\documentclass{amsproc}
\usepackage{amssymb}

\usepackage{graphicx}
\usepackage{amscd}

\theoremstyle{plain}

\newtheorem{corollary}{Corollary}
\newtheorem{definition}{Definition}
\newtheorem{example}{Example}

\newtheorem{proposition}{Proposition}
\newtheorem{remark}{Remark}
\newtheorem{theorem}{Theorem}
\numberwithin{equation}{section}

\begin{document}
\title[Problem on Surjections]{On the S. Banach Problem on Surjections}
\author{Eugene V. Tokarev}
\address{B.E. Ukrecolan, 33-81 Iskrinskaya str., 61005, Kharkiv-5, UKRAINE}
\email{tokarev@univer.kharkov.ua}
\dedicatory{Dedicated to the memory of S. Banach.}

\begin{abstract}
Is shown that any separable superreflexive Banach space $X$ may be
isometricaly embedded in a separable superreflexive Banach space $Z=Z(X)$
(which, in addition, is of the same type and cotype as $X$) such that its
conjugate $Z^{\ast }$ admits a continuous surjection on each its subspace
(and on each separable space $Y^{\ast }$ that is conjugate to some space $Y$%
, which is finitely representable in $X$). This gives an affirmative answer
on S. Banach's problem:

\textit{Whether there exists a Banach space $X$, non isomorphic to a Hilbert
space, which admits a continuous linear surjection on each its subspace and
is essentially different from }$l_{1}$?
\end{abstract}

\maketitle

\section{Introduction}

A long standing problem, which is contained in remarks to the last section
of S.~Banach's monograph [1] is:

\textit{Whether there exists a Banach space $X$, non isomorphic to a Hilbert
space, which admits a continuous linear surjection on each its subspace and
is essentially different from }$l_{1}$?

Such a property of Banach spaces will be referred as to the \textit{%
surjection property.}

Of course, if $X=W\oplus l_{1}$, where $W$ is an arbitrary separable Banach
space, then $X$ enjoys the surjection property. However, up today there are
known no any other examples of Banach spaces with this property.

In the article it will be shown that any separable superreflexive Banach
space $X$ may be isometricaly embedded in a separable superreflexive Banach
space $Z=Z(X)$ (which, in addition, is of the same type and cotype as $X$)
such that its conjugate $Z^{\ast }$ admits a continuous surjection on each
its subspace (and on each separable space $Y^{\ast }$ that is conjugate to
some space $Y$, which is finitely representable in $X$).

An idea of the construction is quite simple. Consider a superreflexive
separable Banach space $X$. It generates a class $X^{<f}$ that consists of
all spaces which are finitely representable in $X$ (exact definitions will
be given below).

Assume that a class $X^{<f}$ is quotient-closed, i.e. that any quotient $Y/Z
$ of any space $Y$, which is finitely representable in $X$, also is finitely
representable in $X$.

Consider a some \textit{envelope} \ $X_{0}$ of this class - i.e. a such
Banach space $X_{0}\in X^{<f}$ that any $Y\in X^{<f}$ of dimension $%
dim(Y)\leq dim(X_{0})$ is isometric to a subspace of $X_{0}$. The existence
of envelopes in every class $X^{<f}$ was proved in additional assumption of
the continuum hypothesis (CH) in [2]; the (CH) from this result was
eliminated in [3].

Immediately, $\left( X_{0}\right) ^{\ast }$ has the surjection property.

Indeed, for any subspace $W\hookrightarrow X$ its conjugate $W^{\ast }$ is
isometric to a quotient of $X_{0}$: $W^{\ast }=\left( X_{0}\right) ^{\ast
\ast }/W^{\perp }=X_{0}/W^{\perp }$ (since $X$ is superreflexive). Since $%
X^{<f}$ is assumed to be quotient closed, $X_{0}/W^{\perp }\in X^{<f}$ and,
since its dimension is not large then that of $X_{0}$ and $X_{0}$ was
assumed to be an envelope, $W^{\ast }$ is isometric to a subspace of $X_{0}$%
. Hence its conjugate $W^{\ast \ast }=W$ is isometric to a quotient of $%
X_{0} $. Since $W\hookrightarrow $ $X_{0}$ was arbitrary, $X_{0}$ has the
surjection property.

So, the main problem is: to built such quotient closed classes $X^{<f}$.

As it will be shown below, for any superreflexive Banach space $Y$ it may be
constructed a such quotient closed class $X^{<f}$ that $Y$ is finitely
representable in $X$; $X$ is of the same type and cotype as $Y$ and, in
addition, class $X^{<f}$ contains a separable space $X_{0}$, which may be
called an \textit{approximate envelope}: for every $\varepsilon >0$ and any
separable space $Z\in X^{<f}$ there exists an isomorphic embedding $%
i:Z\rightarrow X_{0}$ such that $\left\| i\right\| \left\| i^{-1}\right\|
\leq 1+\varepsilon $. As was mention before, a (separable) conjugate space $%
\left( X_{0}\right) ^{\ast }$enjoys the surjection property.

Below will be shown that constructed spaces have nice symmetric (or, rather,
almost symmetric) properties, that are dual to some properties of classical
spaces of almost universal disposition.

\section{Definitions and notations}

\begin{definition}
Let $X$, $Y$ are Banach spaces. $X$ is \textit{finitely representable} in $Y$
(in symbols: $X<_{f}Y$) if for each $\varepsilon >0$ and for every finite
dimensional subspace $A$ of $X$ there exists a subspace $B$ of $Y$ and an
isomorphism $u:A\rightarrow B$ such that $\left\| u\right\| \left\|
u^{-1}\right\| \leq 1+\varepsilon $.

Spaces $X$ and $\ Y$ are said to be finitely equivalent if $X<_{f}Y$ and $%
Y<_{f}X$.

Any Banach space $X$ generates classes
\begin{equation*}
X^{f}=\{Y\in \mathcal{B}:X\sim _{f}Y\}\text{ \ and \ }X^{<f}=\{Y\in \mathcal{%
B}:Y<_{f}X\}
\end{equation*}
\end{definition}

For any two Banach spaces $X$, $Y$ their \textit{Banach-Mazur distance }is
given by
\begin{equation*}
d(X,Y)=\inf \{\left\| u\right\| \left\| u^{-1}\right\| :u:X\rightarrow Y\},
\end{equation*}
where $u$ runs all isomorphisms between $X$ and $Y$ and is assumed, as
usual, that $\inf \varnothing =\infty $.

It is well known that $\log d(X,Y)$ defines a metric on each class of
isomorphic Banach spaces. A set $\frak{M}_{n}$ of all $n$-dimensional Banach
spaces, equipped with this metric, is a compact metric space that is called
\textit{the Minkowski compact} $\frak{M}_{n}$.

A disjoint union $\cup \{\frak{M}_{n}:n<\infty \}=\frak{M}$ is a separable
metric space, which is called the \textit{Minkowski space}.

Consider a Banach space $X$. Let $H\left( X\right) $ be a set of all its
different finite dimensional subspaces (isometric finite dimensional
subspaces of $X$ in $H\left( X\right) $ are identified). Thus, $H\left(
X\right) $ may be regarded as a subset of $\frak{M}$, equipped with the
restriction of the metric topology of $\frak{M}$.

Of course, $H\left( X\right) $ need not to be a closed subset of $\frak{M}$.
Its closure in $\frak{M}$ will be denoted $\overline{H\left( X\right) }$.
From definitions follows that $X<_{f}Y$ if and only if $\overline{H\left(
X\right) }\subseteq \overline{H\left( Y\right) }$. Spaces $X$ and $Y$ are
\textit{finitely equivalent }(in symbols: $X\sim _{f}Y$) if simultaneously $%
X<_{f}Y$ and $Y<_{f}X$. Therefore, $X\sim _{f}Y$ if and only if $\overline{%
H\left( X\right) }=\overline{H\left( Y\right) }$.

There is a one to one correspondence between classes of finite equivalence $%
X^{f}=\{Y\in \mathcal{B}:X\sim _{f}Y\}$ and closed subsets of $\frak{M}$ of
kind $\overline{H\left( X\right) }$.

Indeed, all spaces $Y$ from $X^{f}$ have the same set $\overline{H\left(
X\right) }$. This set, uniquely determined by $X$ (or, equivalently, by $%
X^{f}$), will be denoted by $\frak{M}(X^{f})$ and will be referred as to
\textit{the Minkowski's base of the class} $X^{f}$.

Using this correspondence, it may be defined a set $f\left( \mathcal{B}%
\right) $ of all different classes of finite equivalence, assuming (to
exclude contradictions with the set theory) that members of $f\left(
\mathcal{B}\right) $ are sets $\frak{M}(X^{f})$. For simplicity it may be
says that members of $f\left( \mathcal{B}\right) $ are classes $X^{f}$
itself.

Clearly, $f\left( \mathcal{B}\right) $ is partially ordered by the relation $%
\frak{M}(X^{f})\subseteq \frak{M(}Y^{f})$, which may be replaced by the
relation $X^{f}<_{f}Y^{f}$ of the same meaning. The minimal (with respect to
this order) element of $\ f\left( \mathcal{B}\right) $ is the class $\left(
l_{2}\right) ^{f}$ (the Dvoretzki theorem); the maximal one - the class $%
\left( l_{\infty }\right) ^{f}$ (an easy consequence of the Hahn-Banach
theorem). Other $l_{p}$'s are used in the classifications of Banach spaces,
which was proposed by L. Schwartz [4].

For a Banach space $X$ its $l_{p}$-\textit{spectrum }$S(X)$ is given by
\begin{equation*}
S(X)=\{p\in\lbrack0,\infty]:l_{p}<_{f}X\}.
\end{equation*}

Certainly, if $X\sim_{f}Y$ then $S(X)=S(Y)$. Thus, the $l_{p}$-spectrum $%
S(X) $ may be regarded as a property of the whole class $X^{f}$. So,
notations like $S(X^{f})$ are of obvious meaning.

Let $X$ be a Banach space. It is called:

\begin{itemize}
\item  $c$-\textit{convex,} if $\infty \notin S(X)$;

\item  $B$-\textit{convex,} if $1\notin S\left( X\right) $;

\item  \textit{Finite universal,} if $\infty \in S(X)$.

\item  \textit{Superreflexive,} if every space of the class $X^{f}$ is
reflexive.
\end{itemize}

Equivalently, $X$ is superreflexive if any $Y<_{f}X$ is reflexive. Clearly,
any superreflexive Banach space is $B$-convex.

\section{Quotient closed and divisible classes of finite equivalence}

In this section it will be shown how to enlarge a Minkowski's base $\frak{M}%
(X^{f})$ of a certain $B$-convex class $X^{f}$ to obtain a set $\frak{N}$,
which will be a Minkowski's base $\frak{M}(W^{f})$ for some class $W^{f}$
that holds the $B$-convexity (and, also, $l_{p}$-spectrum) of the
corresponding class $X^{f}$ and will be quotient closed.

\begin{definition}
A class $X^{f}$ (and its Minkowski's base $\frak{M}(X^{f})$) is said to be
divisible if some (equivalently, any) space $Z\in X^{f}$ is finitely
representable in any its subspace of finite codimension.
\end{definition}

\begin{definition}
Let $\{X_{i}:i\in I\}$ be a collection of Banach spaces. A space
\begin{equation*}
l_{2}\left( X_{i},I\right) =\left( \sum \oplus \{X_{i}:i\in I\}\right) _{2}
\end{equation*}
is a Banach space of all families $\{x_{i}\in X_{i}:i\in I\}=\frak{x}$, with
a finite norm
\begin{equation*}
\left\| \frak{x}\right\| _{2}=\sup \{(\sum \{\left\| x_{i}\right\|
_{X_{i}}^{2}:i\in I_{0}\})^{1/2}:I_{0}\subset I;\text{ \ \ }card\left(
I_{0}\right) <\infty \}.
\end{equation*}
\end{definition}

\begin{example}
Any Banach space $X$ may be isometricaly embedded into a space
\begin{equation*}
l_{2}(X)=(\sum\nolimits_{i<\infty }\oplus X_{i})_{2},
\end{equation*}
where all $X_{i}$'s are isometric to $X$. Immediately, $l_{2}(X)$ generates
a divisible class $\mathsf{D}_{2}(X^{f})=\left( l_{2}(X)\right) ^{f}$ which
with the same $l_{p}$-spectrum as $X^{f}$ and is superreflexive if and only
if $X^{f}$ is superreflexive too.
\end{example}

\begin{remark}
The procedure $\mathsf{D}_{2}:X^{f}\rightarrow \left( l_{2}(X)\right) ^{f}$
may be regarded as a closure operator on the partially ordered set $f\left(
\mathcal{B}\right) $. Indeed, it is

\begin{itemize}
\item  Monotone, i.e., $X^{f}<_{f}\mathsf{D}_{2}(X^{f})$;

\item  Idempotent, i.e., $\mathsf{D}_{2}(X^{f})=\mathsf{D}_{2}(\mathsf{D}%
_{2}(X^{f}))$;

\item  Preserve the order: $X^{f}<_{f}Y^{f}$ $\Longrightarrow \mathsf{D}%
_{2}(X^{f})<_{f}\mathsf{D}_{2}(Y^{f})$.
\end{itemize}

It is of interest that extreme points of $f\left( \mathcal{B}\right) $ are
stable under this procedure:
\begin{equation*}
\mathsf{D}_{2}(\left( l_{2}\right) ^{f})=\left( l_{2}\right) ^{f};\text{ }%
\mathsf{D}_{2}(\left( c_{0}\right) ^{f})=\left( c_{0}\right) ^{f}.
\end{equation*}
\end{remark}

To distinguish between general divisible classes and classes of type $%
\mathsf{D}_{2}(X^{f})$, the last ones will be called \textit{2-divisible
classes}.

\begin{definition}
A class $X^{f}$ (and its Minkowski's base $\frak{M}(X^{f})$) is said to be
quotient closed if for any $A\in \frak{M}(X^{f})$ and its subspace $%
B\hookrightarrow A$ the quotient $A/B$ belongs to $\frak{M}(X^{f})$.
\end{definition}

Let $K\subseteq \frak{M}$ be a class of finite dimensional Banach spaces
(recall, that isometric spaces are identified). Define operations $H$, $Q$
and $\ast $ that transform a class $K$ in another class - $H(K)$; $Q(K)$ or $%
(K)^{\ast }$ respectively. Namely, let

\begin{align*}
H(K) & =\{A\in\frak{M}:A\hookrightarrow B;\text{ \ }B\in K\} \\
Q(K) & =\{A\in\frak{M}:A=B/F;\text{ \ }F\hookrightarrow B;\text{ \ }B\in K\}
\\
(K)^{\ast} & =\{A^{\ast}\in\frak{M}:A\in K\}
\end{align*}

In words, $H(K)$ consists of all subspaces of spaces from $K$; $Q(K)$
contains all quotient spaces of spaces of $K$; $(K)^{\ast}$ contains all
conjugates of spaces of $K$.

The following theorem lists properties of these operations. In iteration of
the operations parentheses may be omitted.

Thus, $K^{\ast\ast}\overset{def}{=}\left( (K)^{\ast}\right) ^{\ast}$; $HH(K)%
\overset{def}{=}H(H(K))$ and so on.

\begin{theorem}
Any set $K$ of finite dimensional Banach spaces has the following properties:

\begin{enumerate}
\item  $K^{\ast \ast }=K$; $HH(K)=H(K)$; $QQ(K)=Q(K)$;

\item  $K\subset H(K)$; $K\subset Q(K)$;

\item  If $K_{1}\subset K_{2}$ then $H(K_{1})\subset H(K_{2})$ and $%
Q(K_{1})\subset Q(K_{2})$;

\item  $\left( H(K)\right) ^{\ast }=Q(K^{\ast })$; $\left( Q(K)\right)
^{\ast }=H(K^{\ast })$;

\item  $HQ(HQ(K))=HQ(K)$; $QH(QH(K))=QH(K)$.
\end{enumerate}
\end{theorem}

\begin{proof}
1, 2 and 3 are obvious.

4. If $A\in Q(K)$ then $A=B/E$ for some $B\in K$ and its subspace $E$. Thus,
$A^{\ast }$ is isometric to a subspace of $B^{\ast }.$ Hence, $A^{\ast }\in
H(B^{\ast })$, i.e., $A^{\ast }\in H(K^{\ast })$. Since $A$ is arbitrary, $%
\left( Q(K)\right) ^{\ast }\subseteq H(K^{\ast })$. Analogously, if $B\in K$
and $A\in H(B)$ then $A^{\ast }$ may be identified with a quotient $B^{\ast
}/A^{\perp }$, where $A^{\perp }$ is the annihilator of $A$ in $B^{\ast }$:
\begin{equation*}
A^{\perp }=\{f\in B^{\ast }:\left( \forall a\in A\right) \text{ \ }\left(
f\left( a\right) =0\right) \}.
\end{equation*}

Hence, $A^{\ast }\in (Q(K^{\ast }))^{\ast }$, and thus $\left( H(K)\right)
^{\ast }\subseteq Q(K^{\ast })$.

From the other hand,
\begin{align*}
H(K^{\ast })& =\left( H(K^{\ast })\right) ^{\ast \ast }\subseteq \left(
Q(K^{\ast \ast })\right) ^{\ast }=\left( Q(K)\right) ^{\ast }; \\
Q(K^{\ast })& =\left( Q(K^{\ast })\right) ^{\ast \ast }\subseteq \left(
H(K^{\ast \ast })\right) ^{\ast }=\left( H(K)\right) ^{\ast }.
\end{align*}

5. Let $A\in HQ(K)$. Then $A$ is isometric to a subspace of some quotient
space $E/F$, where $E\in K$; $F\hookrightarrow E$. If $B$ is a subspace of $%
A $ then $\left( A/B\right) ^{\ast }=\left( E/F\right) ^{\ast }/B^{\perp }$,
i.e. $\left( Q\left( HQ(K)\right) \right) ^{\ast }\subseteq Q\left(
Q(K)^{\ast }\right) $. Because of
\begin{align*}
Q\left( HQ(K)\right) & =\left( Q\left( HQ(K)\right) \right) ^{\ast \ast
}\subseteq (Q\left( Q(K)^{\ast }\right) )^{\ast } \\
& \subseteq H(Q(K)^{\ast \ast })=HQ(K),
\end{align*}

we have
\begin{equation*}
H(Q(H(Q(K))))\subseteq H(H(Q(K)))=HQ(K).
\end{equation*}

Analogously, if $A\in QH(K)$, then $A$ is isometric to a quotient space $F/E$%
, where $F\in H(B)$ for some $B\in K$ and $E\hookrightarrow B$. If $W\in
H(A) $, i.e., if $W\in H(F/E)$, then $W^{\ast }=(F/E)^{\ast }/W^{\perp }$
and $(F/E)^{\ast }$ is isometric to a subspace $E^{\perp }$ of $F^{\ast }\in
(H(B))^{\ast }$. Thus, $(H(QH(K)))^{\ast }\subseteq H((H(K))^{\ast })$ and
\begin{align*}
H\left( QH(K)\right) & =\left( H\left( QH(K)\right) \right) ^{\ast \ast
}\subseteq (H\left( H(K)^{\ast }\right) )^{\ast } \\
& \subseteq Q(H(K)^{\ast \ast })=QH(K).
\end{align*}

Hence,
\begin{equation*}
Q(H(Q(HQ(K))))\subseteq Q(Q(H(K)))=QH(K).
\end{equation*}

Converse inclusion follows from 2.
\end{proof}

Consider for a given 2-divisible class $X^{f}$ its Minkowski's base $\frak{M}%
(X^{f})$ and enlarge it by addition to $\frak{M}(X^{f})$ of all quotient
spaces of spaces from $\frak{M}(X^{f})$ and all their subspaces. In the
formal language, consider a set $H(Q(\frak{M}(X^{f})))$.

For any class $W^{f}$ the set $\frak{N}=\frak{M}(W^{f})$ has following
properties:

(\textbf{C}) $\frak{N}\ $is a closed subset of the Minkowski's space $\frak{M%
}$;

(\textbf{H}) If $A\in\frak{N}$ and $B\in H(A)$ then $B\in\frak{N}$;

(\textbf{A}$_{0}$) For any $A$, $B\in\frak{N}$ there exists $C\in\frak{N}$
such that $A\in H(C)$ and $B\in H(C)$.

Clearly, $H(Q(\frak{M}(X^{f})))$ has properties (\textbf{C}), (\textbf{H})
and (\textbf{A}$_{0}$).

\begin{theorem}
Let $\frak{N}$ be a set of finite dimensional Banach spaces; $\frak{N}%
\subset \frak{M}$. If $\frak{N}$ has properties (\textbf{C}), (\textbf{H})
and (\textbf{A}$_{0}$) then there exists a class $X^{f}$ such that $\frak{N}=%
\frak{M}(X^{f})$.
\end{theorem}

\begin{proof}
Since (\textbf{A}$_{0}$), all spaces from $\frak{N}$ may be directed in an
inductive isometric system. Its direct limit $W$ generates a class $W^{f}$
with $\frak{M}(W^{f})\supseteq \frak{N}$. As it follows from the property (%
\textbf{C}), $\frak{M}(W^{f})$ contains no spaces besides that of $\frak{N}$%
. Hence, $\frak{M}(W^{f})=\frak{N}$.

Another proof may be given by using an ultraproduct of all spaces from $%
\frak{N}$ by the ultrafilter, coordinated with a partial ordering on $\frak{N%
}$, which is generated by the property (A$_{0}$).
\end{proof}

Let $X$ be a $B$-convex Banach space.

Let $X=l_{2}(Y)$ (and, hence, $X^{f}=\mathsf{D}_{2}(Y^{f})$). Consider the
Minkowski's base $\frak{M}(X^{f})$ and its enlargement $H(Q(\frak{M}%
(X^{f})))=HQ\frak{M}(X^{f})$.

\begin{theorem}
There exists a Banach space $W$ such that $HQ\frak{M}(X^{f})=\frak{M}(W^{f})$%
.
\end{theorem}

\begin{proof}
Obviously, $HQ\frak{M}(X^{f})$ has properties (H) and (C).

Since $\frak{M}(X^{f})$ is 2-divisible, then for any $A,B\in \frak{M}(X^{f})$
a space $A\oplus _{2}B$ belongs to $\frak{M}(X^{f})$ and, hence, to $HQ\frak{%
M}(X^{f})$.

If $A,B\in Q\frak{M}(X^{f})$ then $A=F/F_{1}$; $B=E/E_{1}$ for some $E,F\in Q%
\frak{M}(X^{f})$.

Clearly, $F/F_{1}\oplus _{2}E$ is isometric to a space ($F\oplus
_{2}E)/F_{1}^{\prime }$, where
\begin{equation*}
F_{1}^{\prime }=\{(f,0)\in F\oplus E:f\in F_{1}\}
\end{equation*}
and, hence, belongs to $Q\frak{M}(X^{f})$. Thus,
\begin{equation*}
F/F_{1}\oplus _{2}E/E_{1}=(F/F_{1}\oplus _{2}E)/E_{1}^{\prime },
\end{equation*}
where
\begin{equation*}
E_{1}^{\prime }=\{(o,e)\in F/F_{1}\oplus _{2}E:e\in E_{1}\},
\end{equation*}
and, hence, belongs to $Q\frak{M}(X^{f})$ too.

If $A,B\in HQ\frak{M}(X^{f})$ then $A\hookrightarrow E$, $B\hookrightarrow F$
for some $E,F\in Q\frak{M}(X^{f})$.

Since $E\oplus _{2}F\in Q\frak{M}(X^{f})$, also $A\oplus _{2}B\in HQ\frak{M}%
(X^{f})$. Thus, $HQ\frak{M}(X^{f})$ has a property (A$_{0}$). A desired
result follows from the preceding theorem.
\end{proof}

\begin{definition}
Let $X$ be a Banach space, which generates a class $X^{f}$\ of finite
equivalence. A class $\ast \ast (X^{f})$ is defined to be a such class $%
W^{f} $ that $\frak{M}(W^{f})=HQ\frak{M}(\mathsf{D}_{2}Y^{f})$
\end{definition}

Clearly, $W^{f}$ is quotient closed. Obviously, $X^{f}<_{f}W^{f}$. It will
be said that $W^{f}$ is a result of a procedure $\ast \ast $ that acts on $f(%
\mathcal{B})$.

\begin{remark}
The procedure $\ast \ast :X^{f}\rightarrow \ast \ast \left( X^{f}\right)
=W^{f}$ may be regarded as a closure operator on the partially ordered set $%
f\left( \mathcal{B}\right) $. Indeed, it is

\begin{itemize}
\item  Monotone, i.e., $X^{f}<_{f}\ast \ast (X^{f})$;

\item  Idempotent, i.e., $X^{f}=\ast \ast (\ast \ast (X^{f}))$;

\item  Preserve the order: $X^{f}<_{f}Y^{f}$ $\Longrightarrow \ast \ast
(X^{f})<_{f}\ast \ast (Y^{f})$.
\end{itemize}

It is of interest that extreme points of $f\left( \mathcal{B}\right) $ are
stable under this procedure:
\begin{equation*}
\ast \ast (\left( l_{2}\right) ^{f})=\left( l_{2}\right) ^{f};\text{\ }\ast
\ast (\left( c_{0}\right) ^{f})=\left( c_{0}\right) ^{f}.
\end{equation*}
\end{remark}

\begin{theorem}
For any Banach space $X$ a class $\ast \ast \left( X^{f}\right) $ is
2-divisible.
\end{theorem}

\begin{proof}
Let $\frak{N}=\frak{M}(\ast \ast (X^{f}))$\ . Since for any pair $A,B\in
\frak{N}$ \ their $l_{2}$-sum belongs to $\frak{N}$, then by an induction, $%
(\sum\nolimits_{i\in I}\oplus A_{i})_{2}\in \frak{N}$ for any finite subset $%
\{A_{i}:i\in I\}\subset \frak{N}$.

Hence, any infinite direct $l_{2}$-sum $(\sum\nolimits_{i\in I}\oplus
A_{i})_{2}$ is finite representable in $\ast \ast \left( X^{f}\right) $. Let
$\{A_{i}:i<\infty \}\subset \frak{N}$ \ is dense in $\frak{N}$. Let $%
Y_{1}=(\sum\nolimits_{i<\infty }\oplus A_{i})_{2}$; $Y_{n+1}=Y_{n}\oplus
_{2}Y_{1}$; $Y_{\infty }=\underset{\rightarrow }{\lim }Y_{n}$. Clearly, $%
Y_{\infty }=l_{2}\left( Y\right) $ and belongs to $\ast \ast \left(
X^{f}\right) $.
\end{proof}

Let $\star :f\left( \mathcal{B}\right) \rightarrow f\left( \mathcal{B}%
\right) $ be one else procedure that will be given by following steps.

Let $X\in \mathcal{B}$; $Y^{f}=\mathsf{D}_{2}\left( X^{f}\right) $. Let $%
\frak{Y}_{0}$ be a countable dense subset of $\frak{M}(Y^{f})$; $\frak{Y}%
_{0}=\left( Y_{n}\right) _{n<\infty }$. Consider a space $Z=(\sum_{n<\infty
}\oplus Y_{n})_{2}$ and its conjugate $Z^{\ast }$. $Z^{\ast }$ generates a
class $\left( Z^{\ast }\right) ^{f}$ which will be regarded as a result of
acting of a procedure $\star :X^{f}\rightarrow \left( Z^{\ast }\right) ^{f}$%
. Since $\left( Z^{\ast }\right) ^{f}$ is 2-divisible, iterations of the
procedure $\star $ are given by following steps: Let $\frak{Z}_{0}$ is a
countable dense subset of $\frak{M}(\left( Z^{\ast }\right) ^{f})$; $\frak{Z}%
_{0}=\left( Z_{n}\right) _{n<\infty }$. Consider a space $W=(\sum_{n<\infty
}\oplus Z_{n})_{2}$ and its conjugate $W^{\ast }$. Clearly,\ $W^{\ast }$
generates a class
\begin{equation*}
\left( W^{\ast }\right) ^{f}=\star \left( Z^{\ast }\right) ^{f}=\star \star
\left( X^{f}\right) .
\end{equation*}

\begin{theorem}
For any Banach space $X$ classes $\ast \ast \left( X^{f}\right) $ and $\star
\star \left( X^{f}\right) $ are identical.
\end{theorem}

\begin{proof}
From the construction follows that $H(Z^{\ast })=\left( QH(l_{2}(X)\right)
^{\ast }$ and that
\begin{align*}
H(W^{\ast })& =(QH(Z^{\ast }))^{\ast }=(Q(QH(l_{2}(X)))^{\ast })^{\ast } \\
& =H(QH(l_{2}(X)))^{\ast \ast }=HQH(l_{2}(X)).
\end{align*}

Hence,
\begin{equation*}
\frak{M}(\left( W^{\ast }\right) ^{f})=HQ\frak{M}(Y^{f})=HQ\frak{M}(\mathsf{D%
}_{2}(X^{f})).
\end{equation*}
\end{proof}

\begin{theorem}
Let $X$ be a Banach space, which generates a class of finite equivalence $%
X^{f}$. If $X$ is B-convex, then the procedure $\ast \ast $ maps $X^{f}$ to
a class $\ast \ast \left( X^{f}\right) $ with the same lower and upper
bounds of its $l_{p}$-spectrum as $X^{f}$. If $X$ is superreflexive then $%
\ast \ast \left( X^{f}\right) $ is superreflexive too.
\end{theorem}

\begin{proof}
Obviously, if $p\in S(X)$ and $p\leq 2$ then the whole interval $[p,2]$
(that may consist of one point) belongs to $S(X)$ and, hence, to $S(\ast
\ast \left( X^{f}\right) )$. If $p\in S(X)$ and $p>2$ then, by duality, $%
p/(p-1)\in S\left( \star \left( X^{f}\right) \right) $; hence $%
[p/(p-1),2]\subset S\left( \star \left( X^{f}\right) \right) $ and, thus $%
[2,p]\subset S(X)\subset S(\ast \ast \left( X^{f}\right) )$. If $p\notin
S(X) $ then $p/(p-1)\notin S\left( \star \left( X^{f}\right) \right) $ by
its construction and, hence, $p\notin S(\ast \ast \left( X^{f}\right) )$ by
the same reason. Hence, if
\begin{equation*}
p_{X}=\inf S(X)\text{; \ }q_{X}=\sup S(X)
\end{equation*}
then
\begin{equation*}
S(\ast \ast \left( X^{f}\right) )=[p_{X},q_{X}]=[\inf S(X),\sup S(X)].
\end{equation*}

The second assertion of the theorem is also obvious.
\end{proof}

\begin{remark}
If $X$ is not $B$-convex, then
\begin{equation*}
\star \left( X^{f}\right) =\star \star \left( X^{f}\right) =\ast \ast \left(
X^{f}\right) =\left( c_{0}\right) ^{f}.
\end{equation*}
\end{remark}

\section{Bifinite representability and properties of classes of kind $\ast
\ast \left( X^{f}\right) $}

To show that classes $\ast \ast \left( X^{f}\right) $ (or, more exactly,
classes $\left( \ast \ast \left( X^{f}\right) \right) ^{<f}$ of all spaces $%
Y $ that are finitely representable in $\ast \ast \left( X^{f}\right) $,
i.e. in any space of $\ast \ast \left( X^{f}\right) $) are indeed quotient
closed in a sense of the introduction it will be needed some more
definitions.

\begin{definition}
(Cf. [9]) Let $X$, $Y$ be Banach spaces. $X$ is said to be bifinitely
representable in $Y$ (shortly, $X<_{ff}Y$) if for every $\varepsilon >0$ and
for every pair of finite dimensional subspaces $A$ of $X$ \ and $A^{\prime }$%
of $X^{\ast }$ there exists a pair of subspaces $B\hookrightarrow Y$ and $%
B^{\prime }\hookrightarrow Y^{\ast }$\ and a pair if isomorphisms $%
u:A\rightarrow B$ and $v:A^{\prime }\rightarrow B^{\prime }$ such that $%
\left\| u\right\| \left\| u^{-1}\right\| \leq 1+\varepsilon $; $\left\|
v\right\| \left\| v^{-1}\right\| \leq 1+\varepsilon $ and $\left\langle
ux,vx^{\prime }\right\rangle =\left\langle x,x^{\prime }\right\rangle $ for
any pair $x\in A$ and $x^{\prime }\in A^{\prime }$.
\end{definition}

To formulate a criterion for validity a relation $X<_{ff}Y$ that was given
in [9] and to prove a desired results will be used a notion on ultrapowers.

\begin{definition}
Let $I$ be a set; $D$ be an ultrafilter over $I$; $\{X_{i}:i\in I\}$ be a
family of Banach spaces. An \textit{ultraproduct } $(X_{i})_{D}$ is given by
a quotient space
\begin{equation*}
(X)_{D}=l_{\infty }\left( X_{i},I\right) /N\left( X_{i},D\right) ,
\end{equation*}
where $l_{\infty }\left( X_{i},I\right) $ is a Banach space of all families $%
\frak{x}=\{x_{i}\in X_{i}:i\in I\}$, for which
\begin{equation*}
\left\| \frak{x}\right\| =\sup \{\left\| x_{i}\right\| _{X_{i}}:i\in
I\}<\infty ;
\end{equation*}
$N\left( X_{i},D\right) $ is a subspace of $l_{\infty }\left( X_{i},I\right)
$, which consists of such $\frak{x}$'s$\ $that
\begin{equation*}
\lim_{D}\left\| x_{i}\right\| _{X_{i}}=0.
\end{equation*}
If all $X_{i}$'s are all equal to a space $X\in \mathcal{B}$ then an
ultraproduct is said to be an \textit{ultrapower} and is denoted by $\left(
X\right) _{D}$.
\end{definition}

In [9] was proved that $X<_{ff}Y$ if there exists a such ultrafilter $D$
that $X^{\ast \ast }$ is isometric to a subspace $Z\hookrightarrow \left(
Y\right) _{D}$ and there exists a projection $P:\left( Y\right)
_{D}\rightarrow Z$ of norm $\left\| P\right\| =1$ (in other words, $X^{\ast
\ast }$ is isometric to an orthogonal complementably subspace of $\left(
Y\right) _{D}$.

\begin{theorem}
Let $X$, $Y$ be Banach spaces. If $X<_{ff}Y$ then $X^{\ast }<_{ff}Y^{\ast }$.
\end{theorem}

\begin{proof}
From the previous discussion follows the existence of a such ultrafilter $D$
that $\left( Y\right) _{D}=X^{\ast \ast }\oplus Z$; $P:\left( Y\right)
_{D}\rightarrow X^{\ast \ast }$; $P^{2}=P$; $\left\| P\right\| =1$. Hence, $%
\left( Y\right) _{D}^{\ast }=X^{\ast \ast \ast }\oplus Z^{\ast }$ and $%
X^{\ast }$ is a range of a projection $P^{\ast }:\left( Y\right) _{D}^{\ast
}\rightarrow X^{\ast }$of norm one. From the principle of local duality of
ultraproducts (cf. [7]) follows that $\left( Y\right) _{D}^{\ast }$ and $%
\left( Y^{\ast }\right) _{D}$ are bifinitely equivalent (i.e. $\left(
Y\right) _{D}^{\ast }<_{ff}$ $\left( Y^{\ast }\right) _{D}$ and $\left(
Y^{\ast }\right) _{D}<_{ff}\left( Y\right) _{D}^{\ast }$; this is
non-trivial fact if $Y$ is not superreflexive: it is known that $\left(
Y\right) _{D}^{\ast }=\left( Y^{\ast }\right) _{D}$ if and only if $Y$ is a
superreflexive Banach space). Thus, $X^{\ast }<_{f}Y^{\ast }$.
\end{proof}

Let $X^{f}$ be a 2-divisible class; $\{B_{i}:i<\infty \}\subset \frak{M}%
(X^{f})$ be dense in $\frak{M}(X^{f})$. Let a space $X_{\oplus }\in X^{f}$
is given by
\begin{equation*}
X_{\oplus }=\left( \sum \oplus \{B_{i}:i<\infty \}\right) _{2}
\end{equation*}

\begin{theorem}
A Banach space $Y$ is finitely representable in $X$ if and only if $Y$ is
bifinite representable in a space $X_{\oplus }$.
\end{theorem}

\begin{proof}
Certainly, $Y<_{ff}X_{\oplus }$ implies that $Y<_{f}X_{\oplus }$ that is
equivalent to $Y<_{f}X$.

To show the converse implication consider a set $H(Y^{\ast \ast })$, which
will be assumed to be indexed by elements of a some set $J$: $H(Y^{\ast \ast
})=\{A_{j}:j\in J\}$.

Let $S_{j}=S_{j}(\varepsilon )=\{i\in J:A_{j}\in H(A_{i})\}\cup
\{\varepsilon \}$. It will be assumed that $S_{j}(\varepsilon _{1})\ll
S_{j}(\varepsilon _{2})$ if $\varepsilon _{1}\geq \varepsilon _{2}$.

For any $n\in \mathbb{N}$ and any choosing $\{i_{1},i_{2},...,i_{n}\}\subset
J$ there exists such $j\in J$ that
\begin{equation*}
S_{j}(\varepsilon )\subset S_{i_{1}}\cap S_{i_{2}}\cap ...\cap S_{i_{n}}.
\end{equation*}

Indeed, from a property (A$_{0}$), which are true for any Minkowski's base
follows that there exists such $j\in J$ that $A_{i_{k}}\in H(A_{j})$ for all
$k=1,2....,n$ .

Therefore a set $\{S_{j}(\varepsilon ):j\in J,$ $\ \varepsilon >0\}$ may be
extended to a some ultrafilter $D$ over $J$, which is concordance with the
order $\ll $.

For $i\in J$ and $\varepsilon >0$ choose a space $B_{n(i)}\in
\{B_{k}:k<\infty \}$ and an isomorphism $u_{i}:A_{i}\rightarrow B_{n\left(
i\right) }$ with $\left\| u_{i}\right\| \left\| u_{i}^{-1}\right\| \leq
1+\varepsilon $. Operators $u_{i}$ may be regarded as embeddings operators $%
u_{i}:A_{i}\rightarrow W$. Let $P_{i}:W\rightarrow B_{n\left( i\right) }$ be
projections of norm one.

Consider an ultrapower $\left( W\right) _{D}$ and define an embedding $%
J:Y^{\ast }\rightarrow \left( W\right) _{D}$ by:
\begin{eqnarray*}
Jy &=&u_{i}^{-1}P_{i}u_{i}\left( y\right) \text{ \ if }y\in A_{i}; \\
&=&0\text{ \ \ \ \ \ \ \ \ \ \ \ \ \ \ \ if }y\notin A_{i}.
\end{eqnarray*}

Let $\left( y_{i}\right) _{i\in I}$ be a family of elements of $W$. Consider
a family $\left( y_{i}^{\prime }\right) _{i\in I}$, which is given by:
\begin{eqnarray*}
y_{i}^{\prime } &=&y_{i}\text{ \ if }y_{i}\in B_{n(i)}; \\
&=&0\text{ \ \ if }y_{i}\notin B_{n(i)}.
\end{eqnarray*}

Let
\begin{equation*}
Q\left( y_{i}\right) _{D}=w^{\ast }-\lim_{D}P_{i}y_{i}^{\prime },
\end{equation*}
where $w^{\ast }-\lim_{D}$ means a limit by an ultrafilter $D$ in a weak*
topology of $Y^{\ast \ast }$. Because of weakly* compactness of the unit
ball of $Y^{\ast \ast }$ this limit exists.

Obviously, the operator $Q\circ J$ is identical on $Y^{\ast \ast }$. The
operator $R=J\circ Q$ defines a projection from $\left( W\right) _{D}$ onto $%
Y^{\ast \ast }$ of norm $\left\| R\right\| =1$. Hence, $Y<_{ff}W$.
\end{proof}

Let $X$ be a Banach space, $X^{f}$ be a corresponding class. $X^{f}$ is said
to be $\ast \ast $-closed if $X^{f}=\ast \ast \left( W^{f}\right) $ for some
$W\in \mathcal{B}$.

\begin{theorem}
If $X^{f}$ is $\ast \ast $-closed, $Y<_{f}X$ and $F\hookrightarrow Y$ then a
quotient $Y/F$ is also finite representable in $X$.
\end{theorem}

\begin{proof}
Consider a dense sequence $\{B_{i}:i<\infty \}\subset \frak{M}(X^{f})$ and a
space
\begin{equation*}
X_{\oplus }=\left( \sum \oplus \{B_{i}:i<\infty \}\right) _{2}.
\end{equation*}

Certainly, $Y<_{ff}X_{\oplus }$. Thus, $Y^{\ast }<_{ff}X_{\oplus }^{\ast }$
, i.e., $Y^{\ast }<_{f}\ast \left( X^{f}\right) $. Of course, $\left(
Y/F\right) ^{\ast }$ is isometric to a subspace of $Y^{\ast }$ and, hence,
is also finite representable in $\ast \left( X^{f}\right) $. Let $%
\{C_{i}:i<\infty \}\subset \frak{M}(\ast X^{f})$ be a dense subsequence. Let
\begin{equation*}
U_{\oplus }=\left( \sum \oplus \{C_{i}:i<\infty \}\right) _{2}.
\end{equation*}
Clearly, $\left( Y/F\right) ^{\ast }<_{ff}U$ and, hence, $\left( Y/F\right)
^{\ast \ast }<_{ff}U^{\ast }$. Since $U^{\ast }\in \ast \ast \left(
X^{f}\right) $, then $\left( Y/F\right) ^{\ast \ast }<_{f}X$ and thus $%
Y/F<_{f}X$.
\end{proof}

At the last of this section, let us describe a set $\frak{M}(\star \left(
X^{f}\right) )$.

\begin{remark}
In the future instead $\star \left( X^{f}\right) $ will be written $\ast
\left( X^{f}\right) $
\end{remark}

\begin{theorem}
A procedure $\ast \left( X^{f}\right) $ asserts to a Minkowski'base $\frak{M}%
(X^{f})$ of a given class $X^{f}$ a set $\frak{M}(\ast \left( X^{f}\right)
)=H\left( \left( \frak{M}(\mathsf{D}_{2}(X^{f}))\right) ^{\ast }\right) $.
\end{theorem}

\begin{proof}
Let $\frak{N}=H\left( \frak{M}(X^{f})\right) ^{\ast }$. If $A,B\in \frak{M}%
(X^{f})$ then $A\oplus _{2}B\in \frak{M}(X^{f})$ and $A^{\ast }\oplus
_{2}B^{\ast }\in \left( \frak{M}(X^{f})\right) ^{\ast }$. If $A,B\in \frak{N}
$ \ then there exists $E$, $F\in \left( \frak{M}(X^{f})\right) ^{\ast }$
such that $A\hookrightarrow E$; $B\hookrightarrow F$. Since $E\oplus
_{2}F\in \frak{N}$, it is obvious that $A\oplus _{2}B\in \frak{N}$ \ too.
Hence, $H\left( \frak{M}(X^{f})\right) ^{\ast }$ satisfies properties (C),
(H) and (A$_{0}$) and is a Minkowski's base of some class $V^{f}$.

Let $A\in \frak{M}(\ast \left( X^{f}\right) )$. Then $A\in H\left(
(U_{\oplus })^{\ast }\right) $, where $U_{\oplus }=\left( \sum \oplus
\{C_{i}:i<\infty \}\right) _{2}$ for a some sequence $\{C_{i}:i<\infty \}$
which is dense in $\frak{M}(X^{f})$. Hence, $A\in H\left( \frak{M}%
(X^{f})\right) ^{\ast }$.

Conversely, if $A\in H\left( \frak{M}(X^{f})\right) ^{\ast }$ then,
obviously, $A\in H\left( (U_{\oplus })^{\ast }\right) \subset \frak{M}(\ast
\left( X^{f}\right) )$.
\end{proof}

\begin{corollary}
For any Banach space $X$ classes $\ast \left( X^{f}\right) $ and $\ast \ast
\ast \left( X^{f}\right) $ are identical.
\end{corollary}

\begin{proof}
Let $\frak{N}=\frak{M}(X^{f})$. Then, according to theorem 1 and its proof,
\begin{eqnarray*}
H\frak{N}^{\ast } &\subseteq &HQ\left( H\frak{N}^{\ast }\right) =HQ\left( Q%
\frak{N}\right) ^{\ast }=H\left( HQ\frak{N}\right) ^{\ast } \\
&=&\left( QHQ\frak{N}\right) ^{\ast }\subseteq \left( HQ\frak{N}\right)
^{\ast }=Q\left( Q\frak{N}\right) ^{\ast }\subseteq H\frak{N}^{\ast }.
\end{eqnarray*}

Since $\frak{M}(\ast \left( X^{f}\right) )=H\frak{N}^{\ast }$, $\frak{M}%
(\ast \ast \left( X^{f}\right) )=HQ\frak{N}$, $\frak{M}(\ast \ast \ast
\left( X^{f}\right) )=H\left( HQ\frak{N}\right) ^{\ast }$, theorem follows.
\end{proof}

\begin{corollary}
Any class $\ast \left( X^{f}\right) $ is $\ast \ast $-closed.
\end{corollary}

\begin{corollary}
In a sequence $X^{f}$; $\ast \left( X^{f}\right) $; $\ast \ast \left(
X^{f}\right) $; $\ast \ast \ast \left( X^{f}\right) $; $\ast \ast \ast \ast
\left( X^{f}\right) $; etc., all even members are identical. Among odd its
members only the first one (namely, $X^{f}$) may differs from others.
\end{corollary}

\begin{remark}
Let $X^{f}$ be a $\ast \ast $-closed class. Classes $X^{f}$ and $\ast \left(
X^{f}\right) $ are in duality in a following sense: $\ast \left( \ast \left(
X^{f}\right) \right) =X^{f}$.
\end{remark}

\section{Spaces of almost universal disposition}

Recall a definition, which is due to V.I. Gurarii [5].

\begin{definition}
Let $X$ be a Banach space; $\mathcal{K}$ be a class of Banach spaces. $X$ is
said to be a space of almost universal disposition with respect to $\mathcal{%
K}$ if for any pair of spaces $A$, $B$ of $\mathcal{K}$ such that $A$ is a
subspace of $B$ ($A\hookrightarrow B$), every $\varepsilon >0$ and every
isomorphic embedding $i:A\rightarrow X$ there exists an isomorphic embedding
$\hat{\imath}:B\rightarrow X$, which extends $i$ (i.e., $\hat{\imath}|_{A}=i$%
) and such, then
\begin{equation*}
\left\| \hat{\imath}\right\| \left\| \hat{\imath}^{-1}\right\| \leq
(1+\varepsilon )\left\| i\right\| \left\| i^{-1}\right\| .
\end{equation*}
\end{definition}

A look on proofs of [5] shows that in the construction of the classical
Gurarii's space of almost universal disposition with respect to $\frak{M}$ \
the main role plays a property of $\frak{M}$, which may be called the
\textit{isomorphic amalgamation property}.

\begin{definition}
Let $X\in \mathcal{B}$ generates a class $X^{f}$\ with a Minkowski's base $%
\frak{M}(X^{f})$. It will be said that $\frak{M}(X^{f})$ (and the whole
class $X^{f}$) has the isomorphic amalgamation property if for any fifth $%
\left\langle A,B_{1},B_{2},i_{1},i_{2}\right\rangle $, where $A$, $B_{1}$, $%
B_{2}\in \frak{M}(X^{f})$; $i_{1}:A\rightarrow B_{1}$ and$\
i_{2}:A\rightarrow B_{2}$ are isomorphic embeddings, there exists a triple $%
\left\langle j_{1},j_{2},F\right\rangle $, where $F\in \frak{M}(X^{f})$ and$%
\ j_{1}:B_{1}\rightarrow F$; $j_{2}:B_{2}\rightarrow F$ are isometric
embedding, such that $j_{1}\circ i_{1}=j_{2}\circ i_{2}$.
\end{definition}

\begin{theorem}
Any 2-divisible quotient closed class $X^{f}$ has the isomorphic
amalgamation property.
\end{theorem}

\begin{proof}
Let $A$, $B_{1}$, $B_{2}\in \frak{M}(X^{f})$; $i_{1}:A\rightarrow B_{1}$ and
$i_{2}:A\rightarrow B_{2}$ are isomorphic embeddings. Since $X^{f}$ is
2-divisible, $C=B_{1}\oplus _{2}B_{2}\in \frak{M}(X^{f})$. Consider a
subspace $H$ of $C$ that is formed by elements of kind $\left(
i_{1}a,-i_{2}a\right) $, where $a$ runs $A$:
\begin{equation*}
H=\{\left( i_{1}a,-i_{2}a\right) :a\in A\}.
\end{equation*}

Consider a quotient $W=C/H$. Since $X^{f}$ is quotient closed, $W\in \frak{M}%
(X^{f})$. Let $h:C\rightarrow W$ be a standard quotient map. Let $%
j_{1}:B_{1}\rightarrow W$ and $j_{2}:B_{2}\rightarrow W$ are given by:
\begin{align*}
j_{1}(b_{1})& =h(b_{1},0)\text{ \ \ \ for }b_{1}\in B_{1}; \\
j_{2}(b_{2})& =h(0,b_{2})\text{ \ \ \ for }b_{2}\in B_{2}.
\end{align*}

It is clear that $j_{1}$ and $j_{2}$ are isometric embeddings such that $\
j_{1}\circ i_{1}=j_{2}\circ i_{2}$.
\end{proof}

\begin{remark}
It may be shown that any quotient closed divisible (not necessary
2-divisible) class $X^{f}$ enjoys the isomorphic amalgamation property too.
For aims of the article the preceding result is satisfactory.
\end{remark}

The proof of a following results is almost literally repeats the Gurarii's
proof [5] on existence of a space of almost universal disposition with
respect to $\frak{M}$. Only changes that need to be made are: a substitution
of a set $\frak{M}$ with $\frak{M}(X^{f})$ for a given class $X^{f}$ and
using instead the mentioned above isomorphic amalgamation property of a set $%
\frak{M}$ \ the same property of a set $\frak{M}(X^{f})$. For this reason
the proof of it is omitted.

\begin{theorem}
Any $\ast \ast $-closed class $X^{f}$ contains a separable space $E_{X}$ of
almost universal disposition with respect to a set $\frak{M}(X^{f})$. This
space is unique up to almost isometry and is almost isotropic (in an
equivalent terminology, has an almost transitive norm: for any two elements $%
a$, $b\in E_{X}$, such that $\left\| a\right\| =\left\| b\right\| $ \ and
every $\varepsilon >0$ there exists an automorphism $u=u(a,b,\varepsilon
):E_{X}\overset{onto}{\rightarrow }E_{X}$ such that $\left\| u\right\|
\left\| u^{-1}\right\| \leq 1+\varepsilon $ and $ua=b$). This space is an
approximative envelope of a class $X^{f}$: for any $\varepsilon >0$ every
separable Banach space which is finitely representable in $X^{f}$ is $%
(1+\varepsilon )$-isomorphic to a subspace of $E_{X}$.
\end{theorem}

\begin{definition}
Let $X\in \mathcal{B}$ generates a class $X^{f}$\ with a Minkowski's base $%
\frak{M}(X^{f})$. It will be said that $\frak{M}(X^{f})$ (and the whole
class $X^{f}$) has the coamalgamation property if for any $\left\langle
A,B_{1},B_{2},h_{1},h_{2}\right\rangle $, where $A$, $B_{1}$, $B_{2}\in
\frak{M}(X^{f})$; $h_{1}:B_{1}\rightarrow A$ and$\ h_{2}:B_{2}\rightarrow A$
are surjections of norm one, there exists a triple $\left\langle
H_{1},H_{2},F\right\rangle $, where $F\in \frak{M}(X^{f})$; $\
H_{1}:F\rightarrow B_{1}$ and $H_{2}:F\rightarrow B_{2}$ are surjections of
norm one such that $h_{1}\circ H_{1}=h_{2}\circ H_{2}$.
\end{definition}

\begin{theorem}
For a $\ast \ast $-closed class $X^{f}$ the dual class $\ast \left(
X^{f}\right) $ has the coamalgamation property.
\end{theorem}

\begin{proof}
Immediately follows from duality and results of the previous section.
\end{proof}

\begin{corollary}
Any $\ast \ast $-closed class $X^{f}$ has the coamalgamation property.
\end{corollary}

\begin{proof}
As it follows from preceding results, any $\ast \ast $-closed class $X^{f}$
is the dual class $\ast W^{f}$, where $W^{f}=\ast \left( X^{f}\right) $, and
hence has the coamalgamation property.
\end{proof}

\begin{definition}
A Banach space $X$ is said to be couniversal for a class $\mathcal{K}$ of
Banach spaces if there exists a continuous linear surjection $h:X\rightarrow
F$ for any $F\in \mathcal{K}$ .
\end{definition}

\begin{theorem}
Let $X^{f}$ be a $\ast \ast $-closed superreflexive class; $E_{X}$ be a
separable space of almost universal disposition with respect to a set $\frak{%
M}(X^{f})$. Then its conjugate, $\left( E_{X}\right) ^{\ast }$ is
couniversal for all separable spaces that are finite representable in $\ast
\left( X^{f}\right) $.
\end{theorem}

\begin{proof}
Let $W<_{f}\ast \left( X^{f}\right) $. Then $W^{\ast }$ is isomorphic (up to
any $\varepsilon >0$) to a subspace of $E_{X}$. Hence, its conjugate $%
W^{\ast \ast }=W$ (by superreflexivity) is isomorphic (up to the same $%
\varepsilon >0$)\ to a quotient $\left( E_{X}\right) ^{\ast }/\left( W^{\ast
}\right) ^{\perp }$.
\end{proof}

\begin{remark}
In a general case $\left( E_{X}\right) ^{\ast }$ is not the unique space
which is couniversal for all separable spaces that are finite representable
in $\ast \left( X^{f}\right) $. Other such spaces may be obtained as
conjugates to other separable approximative envelope of a class $X^{f}$.
\end{remark}

\begin{theorem}
Let $X^{f}$ be a $\ast \ast $-closed class. If all its separable
approximative envelopes are pairwice isomorphic then $X$ is isomorphic to a
Hilbert space.
\end{theorem}

\begin{proof}
Since $E_{X}\oplus _{2}X_{\oplus }$ is also a separable approximative
envelope of a class $X^{f}$, it is sufficient to show that if $E_{X}$ and $%
E_{X}\oplus _{2}X_{\oplus }$ are isomorphic, than both of them are
isomorphic to a Hilbert space.

Let $A\hookrightarrow E_{X}$ be a finite dimensional subspace. Let
\begin{equation*}
\lambda (A\hookrightarrow E_{X})=\inf \{\left\| P\right\| :\text{ \ \ }%
P:E_{X}\rightarrow A\},
\end{equation*}
where $P$ runs all projections of $E_{X}$ onto $A$.

It easy to see that $\lambda (iA\hookrightarrow E_{X})=\lambda
(jA\hookrightarrow E_{X})$ for any pair of isometric embeddings, $i,j$ of $A$
into $E_{X}$. Indeed, as was mentioned before, for any $\varepsilon >0$
there is an automorphism $u$ of $E_{X}$ with $\left\| u\right\| \left\|
u^{-1}\right\| \leq 1+\varepsilon $ such that $i\circ u=j$. Let $%
P:E_{X}\rightarrow A$ be a projection. Then $P\circ u:E_{X}\rightarrow A$is
a projection of norm $\left\| P\circ u\right\| \leq \left\| P\right\|
\left\| u\right\| \left\| u^{-1}\right\| \leq \left( 1+\varepsilon \right)
\left\| P\right\| $. Since $P$ and $\varepsilon $ are arbitrary, then
\begin{equation*}
\lambda (iA\hookrightarrow E_{X})=\lambda (jA\hookrightarrow E_{X}).
\end{equation*}

Consider a space $E_{X}\oplus _{2}X_{\oplus }$. For any $A\in H\left(
X_{\oplus }\right) $%
\begin{equation*}
\lambda (A,E_{X}\oplus X_{\oplus })=\inf \{\lambda (iA\hookrightarrow
E_{X}\oplus _{2}X_{\oplus }):i:A\rightarrow E_{X}\oplus _{2}X_{\oplus }\}=1.
\end{equation*}

If $d\left( E_{X}\oplus X_{\oplus },E_{X}\right) \leq C<\infty $ then, since
$H\left( X_{\oplus }\right) $ is dense in $\frak{M}\left( X^{f}\right) $,
for any finite dimensional subspace $A\hookrightarrow E_{X}$
\begin{equation*}
\lambda ^{\prime }(E_{X})=\sup \{\lambda (iA\hookrightarrow E_{X}):\text{ \ }%
i:A\rightarrow E_{X}\text{; \ }A\in H\left( E_{X}\right) \}\leq C.
\end{equation*}

According to [6], $E_{X}$ is isomorphic to the Hilbert space $l_{2}$.
\end{proof}

\begin{corollary}
Any class $\ast \ast \left( X^{f}\right) $ either is generated be a space
isomorphic to the Hilbert one or contains a continuum number of pairwise
non-isomorphic separable approximative envelopes.
\end{corollary}

\begin{proof}
Let $E_{X}$ be not isomorphic to the Hilbert space $l_{2}$. Then there
exists an infinite sequence $\left( A_{m}\right) _{m<\infty }$ of finite
dimensional subspaces of $E_{X}$, such that $\lambda (A_{m}\hookrightarrow
E_{X})\geq m$ ($m\in \mathbb{N}$). Let $M\in \mathbb{N}$ be an infinite
subset. Consider a space
\begin{equation*}
X_{M}=\sum \oplus \{A_{m}:m\in M\}_{2}.
\end{equation*}

Among spaces of kind $X_{M}$ there are continuum of pairwice non-isomorphic
ones. Any space $X_{M}$ is separable; any sum $E_{X}\oplus _{2}X_{M}$ is an
approximative envelope of $\ast \ast \left( X^{f}\right) $; if $X_{M_{1}}$
and $X_{M_{2}}$ are non-isomorphic then $E_{X}\oplus _{2}X_{M_{1}}$ and $%
E_{X}\oplus _{2}X_{M_{1}}$ are also non isomorphic. This proves the
corollary.
\end{proof}

\begin{theorem}
\ For any $\ast \ast $-closed superreflexive class $X^{f}$ and for the
corresponding separable space $E_{X}$ of almost universal disposition with
respect to a set $\frak{M}(X^{f})$, its conjugate, $\left( E_{X}\right)
^{\ast }$, has the surjection property, i.e. admits a continuous (linear)
surjection on each its subspace.
\end{theorem}

\begin{proof}
$\left( E_{X}\right) ^{\ast }$ and , hence, each its subspace is finite
representable in $\ast \left( X^{f}\right) $. Since $\ast \left(
X^{f}\right) $ is superreflexive and $E_{X}$ is separable then $\left(
E_{X}\right) ^{\ast }$ is separable too. The result follows from the
preceding theorem.
\end{proof}

This theorem solves the mentioned above Banach's problem on surjections.

\begin{remark}
It is of interest, that constructed above spaces $E_{X}$ and their
conjugates $\left( E_{X}\right) ^{\ast }$ has a lot of additional properties
that will be described in other article. Here will be given simplest of them.
\end{remark}

\begin{proposition}
Let $\left( E_{X}\right) ^{\ast }$ be a conjugate to a space of almost
universal disposition from a $\ast \ast $-closed class $X^{f}$ that is
B-convex but not superreflexive. Then $\left( E_{X}\right) ^{\ast }$ does
not contain any uncountable unconditional sequence.
\end{proposition}

\begin{proof}
Let $X^{f}$ be non-superreflexive but B-convex. Assume that $\{x_{i}:i\in
I\}\subset \left( E_{X}\right) ^{\ast }$ be an uncountable unconditional
sequence. Then it spans a reflexive subspace
\begin{equation*}
W=span\{x_{i}:i\in I\}
\end{equation*}
of $\left( E_{X}\right) ^{\ast }$. Clearly, $W$ is weakly* closed and,
hence, there exists a quotient $E_{X}/Z$, which conjugate $\left(
E_{X}/Z\right) ^{\ast }$ is isometric to $W$. However $E_{X}/Z$ is reflexive
and separable, and hence, $\left( E_{X}/Z\right) ^{\ast }$ also has these
properties. This contradicts with non-separability of $W$.
\end{proof}

\begin{remark}
Analogously, $\left( E_{X}\right) ^{\ast }$ contains no non-separable
reflexive subspaces.
\end{remark}

One more property of $\left( E_{X}\right) ^{\ast }$ in a superreflexive case
is the following property, that may be called the \textit{almost
cohomoheneity}.

\begin{definition}
A Banach space $X$ is said to be almost cohomoheneous if for every $%
\varepsilon >0$ and for any pair of finite dimensional isometric quotients $%
X/Y$ and $X/Z$ the corresponding isometry $i:X/Y\rightarrow X/Z$ may be
lifted to an $(1+\varepsilon )$-automorphism of $X$, i.e. there exists such
automorphism $I:X\rightarrow X$ that $\left\| I\right\| \left\|
I^{-1}\right\| \leq 1+\varepsilon $ and that $i\circ h_{Y}=h_{Z}\circ I$,
where $h_{Y}:X\rightarrow X/Y$; $h_{Z}:X\rightarrow X/Z$ are quotient
mappings.
\end{definition}

\begin{proposition}
Let $\left( E_{X}\right) ^{\ast }$ be a conjugate to a space of almost
universal disposition from a $\ast \ast $-closed class $X^{f}$. Then $\left(
E_{X}\right) ^{\ast }$ is almost cohomoheneous.
\end{proposition}

\begin{proof}
Easy follows from definition of $E_{X}$ and from duality.
\end{proof}

\begin{remark}
The existence of a such separable space $W_{X}$ (which is finitely
representable in a superreflexive $\ast \ast $-closed class $X^{f}$) that
posses such properties as couniversality and almost cohomoheneity may be
proved directly, without appealing to $E_{X}$, by using the coamalgamation
property of a set $\frak{M}(X^{f})$.
\end{remark}

It is naturally to called such a space to be a space of \textit{almost
couniversal disposition.}

\section{References}

\begin{enumerate}
\item  Banach S. Th\'{e}orie des op\'{e}rations lin\'{e}aires, Warszawa -
Lw\'{o}w. Z Subwencji funduszu kultury Narodowej. Monografie Matematyczne,
\textbf{1} (1932)

\item  Stern J. The problem of envelopes for Banach spaces, Israel J. Math.
\textbf{24:1} (1976) 1 - 15

\item  Tokarev E.V.\textit{Problem of envelopes of locally convex spaces }%
(transl. from Russian), Siberian Mathematical Journal, \textbf{31:1 }(1990)
173-175

\item  Schwartz L. \textit{Geometry and probability in Banach spaces}, Bull.
AMS \textbf{4:2} (1981) 135-141

\item  Gurarii V.I.\textit{\ Spaces of universal disposition, isotropic
spaces and the Mazur problem on rotations in Banach spaces}, Sibirsk. Mat.
Journ. (in Russian) \textbf{7} (1966) 1002-1013

\item  Lindenstrauss J., Tzafriri L. \textit{On the complemented subspaces
problem}, Israel J. Math. \textbf{9} (1971) 263 - 269
\end{enumerate}

\end{document}